\documentclass[12pt]{amsart}
\usepackage{amscd}     
\usepackage{amssymb}
\usepackage{amsmath, amsthm, graphics}
\usepackage{xypic}     
\LaTeXdiagrams        
\usepackage[all]{xy}
\xyoption{2cell} \UseAllTwocells \xyoption{frame} \CompileMatrices
\allowdisplaybreaks[3]
\usepackage{amsfonts}
\usepackage[normalem]{ulem}
\usepackage{hyperref}
\usepackage{tikz}
\usepackage{mathtools}
\usepackage{verbatim} 

\usepackage[margin=4cm]{geometry}

\usepackage{latexsym}
\usepackage{epsfig}
\usepackage{amsfonts}
\usepackage{enumerate}
\usepackage{times,float}
\usepackage{graphics}

\newtheorem{theorem}{Theorem}[section]

\newtheorem{problem}[theorem]{Problem}

\theoremstyle{remark}

\theoremstyle{remark}

\newtheorem{remark}[subsection]{Remark}

\numberwithin{equation}{section}

\def\b1{{\mathbf 1}}

\begin{document}

\title{A note on Virasoro constraints for products}
\author{Hsian-Hua Tseng}
\address{Department of Mathematics\\ The Ohio State University\\ 100 Math Tower, 231 West 18th Avenue\\ Columbus, OH 43210\\ USA}
\email{hhtseng@math.ohio-state.edu}
\date{\today}

\begin{abstract}
We study Virasoro constraints for Gromov-Witten theory of a product variety when one factor has semi-simple quantum cohomology.
\end{abstract}

\maketitle

\section{Introduction}\label{sec:intro}
Let $W$ be a smooth projective variety over $\mathbb{C}$. The Virasoro conjecture for $W$ is an infinite collection of differential equations for the total descendant potential $\mathcal{D}^W$ of Gromov-Witten theory of $W$:
\begin{equation*}
L_m^W\mathcal{D}^W=0, \quad m\geq -1.    
\end{equation*}
A brief discussion of the differential operators $L_m^W$ is given in Section \ref{sec:vir}.

Virasoro conjecture is an outstanding open problem in Gromov-Witten theory dated back to the 1990s. A review of Virasoro conjecture from the early days can be found in \cite{ge}. For a more recent survey, see \cite{cgt}.

To the best of our knowledge, proven cases of Virasoro conjecture are based on one of the two approaches. The approach of \cite{OP} establishes Virasoro conjecture for nonsingular curves. The approach of Givental \cite{g2} makes use of the structure of {\em semi-simple} Gromov-Witten theory.

Let $X$ and $Y$ be smooth projective varieties over $\mathbb{C}$. The following problem is natural.
\begin{problem}\label{problem_product}
Show that Virasoro conjecture holds for $X\times Y$ if and only if Virasoro conjecture holds for both $X$ and $Y$.    
\end{problem}
Although Problem \ref{problem_product} is natural, there has been very little progress on Problem \ref{problem_product} in the past decades\footnote{The result on toric bundles \cite{cgt} contains some cases of Problem \ref{problem_product}.}. Hoping to create some momentum for progress on Problem \ref{problem_product}, the goal of this note is to present the following partial result on Problem \ref{problem_product}:

\begin{theorem}\label{thm:main}
Suppose the quantum cohomology of $Y$ is semi-simple at some point. Then Virasoro constraints for $X\times Y$ hold if and only if Virasoro constraints hold for $X$.     
\end{theorem}

Theorem \ref{thm:main} is derived by applying Givental's approach: more precisely, the approach to study Virasoro constraints for toric bundles \cite{cgt}. The semi-simplicity assumption on $Y$ allows us to apply Givental-Teleman classification to the Gromov-Witten theory of $Y$. Together with product formula of Gromov-Witten invariants \cite{b}, we obtain a formula for Gromov-Witten theory of $X\times Y$ in terms of Gromov-Witten theory of $X$ and the action of certain loop group element: we first do this at the level of {\em cohomological field theories}, then at the level of generating functions. Once such a formula is obtained, we check that the loop group element involved preserves grading. Theorem \ref{thm:main} then follows from loop group covariance \cite{cgt}.

The rest of this note is organized as follows. In Section \ref{sec:prelim} we summarize the required materials. Section \ref{sec:CohFT} contains a discussion of cohomological field theory. Section \ref{sec:action} contains a summary of the $R$-matrix action (also known as actions by loop group elements). Section \ref{sec:classification} presents a brief summary of the Givental-Teleman classification of semi-simple cohomological field theories. The main result is treated in Section \ref{sec:product_target}. In Section \ref{sec:GW_CohFT}, we review the construction of cohomological field theories from Gromov-Witten theory. In Section \ref{sec:CohFT_product}, we derive formulas for the Gromov-Witten theory of a product $X\times Y$. A brief summary of Virasoro constraints is given in Section \ref{sec:vir}. Finally, Theorem \ref{thm:main} is derived in Section \ref{sec:pf_main}.

\subsection{Acknowledgment}
The author is supported in part by a Simons Foundation collaboration grant.

\section{Preliminary materials}\label{sec:prelim}
\subsection{Cohomological field theory}\label{sec:CohFT}
We recall the notion of cohomological field theories. Our presentation follows closely that of \cite[Section 0.5]{PPZ}, see also \cite{P}.

Let $V$ be a finite dimensional vector space over a field of characteristic $0$, equipped with a non-degenerate pairing $\eta$ and a distinguished element $1$. A {\em cohomological field theory (CohFT)} with unit (modelled on $(V, \eta,1)$) is a system $\Omega=(\Omega_{g,n})_{2g-2+n>0}$ where 
\begin{equation*}
\Omega_{g,n}\in H^*(\overline{\mathcal{M}}_{g,n})\otimes (V^*)^{\otimes n},    
\end{equation*}
subject to the following axioms:
\begin{enumerate}
\item[(i)] Each $\Omega_{g,n}$ is invariant under the action of the symmetric group $S_n$ given by permuting the marked points of $\overline{\mathcal{M}}_{g,n}$ and copies of $V^*$.

\item[(ii)] The pull-backs $q^*\Omega_{g,n}$ and $r^*\Omega_{g,n}$ under the gluing maps
\begin{equation*}
q:\overline{\mathcal{M}}_{g-1,n+2}\to \overline{\mathcal{M}}_{g,n}, \quad r: \overline{\mathcal{M}}_{g_1,n_1+1}\times \overline{\mathcal{M}}_{g_2,n_2+1}\to \overline{\mathcal{M}}_{g,n} 
\end{equation*}
are equal to the contractions of $\Omega_{g-1,n+2}$ and $\Omega_{g_1, n_1+1}\otimes \Omega_{g_2,n_2+1}$ by the bi-vector
\begin{equation*}
\sum_{j,k}\eta^{jk}e_j\otimes e_k    
\end{equation*}
inserted at the two identified points. Here $\{e_i\}\subset V$ is a basis, $\eta_{jk}=\eta(e_j,e_k)$, and $\eta^{jk}$ are matrix coefficients of the inverse of the matrix $(\eta_{jk})$. 

\item[(iii)]
Under the forgetful map $p: \overline{\mathcal{M}}_{g,n+1}\to \overline{\mathcal{M}}_{g,n}$, we have
\begin{equation*}
\Omega_{g,n+1}(v_1\otimes...\otimes v_n\otimes 1)=p^*\Omega_{g,n}(v_1\otimes...\otimes v_n), \quad v_1,...,v_n\in V.    
\end{equation*}
Also, $\Omega_{0,3}(v_1\otimes v_2\otimes 1)=\eta(v_1,v_2)$.
\end{enumerate}

A {\em CohFT} consists of the above data without $1\in V$ and axiom (iii).

A CohFT defines a {\em quantum product} $\bullet$ on $V$ by $\eta(v_1\bullet v_2, v_3)=\Omega_{0,3}(v_1\otimes v_2\otimes v_3).$

A {\em topological field theory} $\omega$ is a CohFT composed only of degree $0$ classes:
\begin{equation*}
\omega_{g,n}\in H^0(\overline{\mathcal{M}}_{g,n})\otimes (V^*)^{\otimes n}.
\end{equation*}

Associated to a CohFT as above, we can define its potential $\mathcal{A}^\Omega$ as follows (see \cite[Equation (6.1)]{t}). For $v(z)=\sum_k v_k z^k\in V[[z]]$, define
\begin{equation}
\mathcal{A}^\Omega(v)=\exp\left(\sum_{g,n}\frac{\hbar^{g-1}}{n!}\int_{\overline{\mathcal{M}}_{g,n}}\Omega_{g,n}(v(\psi_1)\otimes...\otimes v(\psi_n) \right).    
\end{equation}
Here, inserting a class $v(\psi_i)$ into the CohFT $\Omega$ means the following: $\Omega_{g,n}(...v(\psi_i)...)=\sum_{k\geq 2}\phi_i^k\Omega_{g,n}(...v_k...)$.

\subsection{The $R$-matrix action}\label{sec:action}
We recall the essence of the $R$-matrix action on CohFTs following \cite[Section 2]{PPZ}.

Let $\Omega$ be a cohomological field theory (CohFT) with unit modelled on $(V, \eta,1)$, as in Section \ref{sec:CohFT}. Consider the group of $\text{End}(V)$-valued (formal) power series 
\begin{equation*}
R(z)=\text{Id}+R_1z+R_2z^2+...    
\end{equation*}
satisfying the symplectic condition 
\begin{equation*}
R(z)R^*(-z)=\text{Id},    
\end{equation*}
where $R^*$ is the adjoint with respect to $\eta$. 

The expression $R^{-1}(z)=\text{Id}/R(z)$ is defined as a formal power series. The symplectic condition implies that $R^{-1}(z)=R^*(-z)$.

Given such a $R$, we define a CohFT $R\Omega$ as follows. 

Recall (e.g. \cite[Section 0.2]{PPZ}) that a stable graph is the collection 
\begin{equation*}
\Gamma=(\mathrm{V}, \mathrm{H}, \mathrm{L}, \mathrm{g}:\mathrm{V}\to \mathbb{Z}_{\geq 0}, v:\mathrm{H}\to\mathrm{V}, \iota:\mathrm{H}\to\mathrm{H})    
\end{equation*}
satisfying the following properties:
\begin{enumerate}

\item $\mathrm{V}$ is the vertex set with the genus function $\mathrm{g}:\mathrm{V}\to \mathbb{Z}_{\geq 0}$,

\item $\mathrm{H}$ is the set of half-edges with the vertex assignment $v:\mathrm{H}\to\mathrm{V}$ and an involution $\iota:\mathrm{H}\to \mathrm{H}$,

\item $\mathrm{E}$ is the set of edges, which are defined to be $2$-cycles of $\iota$ in $\mathrm{H}$ (including self-edges at vertices),

\item $\mathrm{L}$ is the set of legs, which are defined to be fixed points of $\iota$ and equipped with a bijection with the set of markings,

\item the pair $(\mathrm{V}, \mathrm{E})$ defines a connected graph,

\item for each vertex $v\in \mathrm{V}$, define $\mathrm{n}(v)$ to be the number of edges and legs incident at $v$. The stability condition $2\mathrm{g}(v)-2+\mathrm{n}(v)>0$ is required.

\end{enumerate}

The genus of a stable graph $\Gamma$ is defined by $\mathrm{g}(\Gamma)=\sum_{v\in \mathrm{V}}\mathrm{g}(v)+h^1(\Gamma)$.

For a stable graph $\Gamma$ of genus $g$ with $n$ legs, define
\begin{equation*}
\text{Cont}_\Gamma\in H^*(\overline{\mathcal{M}}_{g,n})\otimes (V^*)^{\otimes n}
\end{equation*}
by the following rules\footnote{We refer to \cite[Section 0.2]{PPZ} for the $\psi$ classes assigned to legs and edges.}:
\begin{enumerate}
    \item place $\Omega_{\mathrm{g}(v),\mathrm{n}(v)}$ at each vertex $v$ of $\Gamma$,

    \item place $R^{-1}(\psi_l)$ at every leg $l$ of $\Gamma$,

    \item at every edge $e$ of $\Gamma$, place 
    \begin{equation*}
    \frac{\eta^{-1}-R^{-1}(\psi_e')\eta^{-1}R^{-1}(\psi_e'')^{\mathbf{t}}}{\psi_e'+\psi_e''}.    
    \end{equation*}
    
\end{enumerate}

Let $\mathsf{G}_{g,n}$ be the (finite) set of stable graphs $\Gamma$ of genus $g$ with $n$ legs. The CohFT $R\Omega$ is defined by 
\begin{equation}
(R\Omega)_{g,n}=\sum_{\Gamma\in\mathsf{G}_{g,n}}\frac{1}{|\text{Aut}(\Gamma)|}\text{Cont}_\Gamma,    
\end{equation}
see \cite[Definition 2.2]{PPZ}.

For a $V$-valued power series without terms in degrees $0$ and $1$,
\begin{equation*}
T(z)=T_2z^2+T_3z^3+...,    
\end{equation*}
define the CohFT $T\Omega$, the translation of $\Omega$ by $T$, by
\begin{equation*}
(T\Omega)_{g,n}(v_1\otimes...\otimes v_n)=\sum_{m\geq 0}\frac{1}{m!}(p_m)_*\Omega_{g,n+m}(v_1\otimes...\otimes v_n\otimes T(\psi_{n+1})\otimes...\otimes T(\psi_{n+m})),    
\end{equation*}
where $p_m:\overline{\mathcal{M}}_{g,n+m}\to \overline{\mathcal{M}}_{g,n}$ is the forgetful map. See \cite[Definition 2.5]{PPZ}.

Suppose $\Omega$ is a CohFT {\em with unit}. For an $R$-matrix $R(z)$ as above, define
\begin{equation*}
T(z)=z\cdot1-zR^{-1}(z)(1)\in z^2 V[[z]].    
\end{equation*}
The {\em unit-preserving $R$-matrix action} on $\Omega$ is defined to be 
\begin{equation*}
R.\Omega=RT\Omega,    
\end{equation*}
see \cite[Definition 2.13]{PPZ}. The point here is that this combination of the $R$-matrix action and translation yields a CohFT with the same unit $1\in V$.

\subsection{Givental-Teleman classification}\label{sec:classification}
Let $\Omega$ be a cohomological field theory with unit modelled on $(V, \eta,1)$. Suppose the algebra $(V,\bullet, 1)$ given by the quantum product of $\Omega$ is semi-simple. Let $\omega$ be the topological field theory associated to $\Omega$ by taking degree $0$ part. The Givental-Teleman classification \cite{g, g2}, \cite{t} may be stated as follows:
\begin{theorem}\label{thm:gt}
There exists a unique 
\begin{equation*}
R(z)=\text{Id}+R_1z+R_2z^2+...\in \text{End}(V)[[z]]    
\end{equation*}
satisfying the symplectic condition $R(z)R^*(-z)=\text{Id}$, such that 
\begin{equation*}
\Omega=R.\omega.    
\end{equation*}
\end{theorem}
See \cite[Section 1.3]{P} for more details.

In the semi-simple case, the topological part $\omega$ can be evaluated explicitly using CohFT axioms. Let $\{e_i\}$ be the idempotent basis of $(V,\bullet, 1)$:
\begin{equation*}
e_i\bullet e_j=\delta_{ij}e_i.    
\end{equation*}
Let 
\begin{equation*}
\tilde{e}_i=e_i/\eta(e_i,e_i)^{1/2}    
\end{equation*}
be the normalized idempotents. We have
\begin{equation}\label{eqn:top_FT}
\omega_{g,n}(\tilde{e}_{i_1}\otimes...\otimes \tilde{e}_{i_n})=
\begin{cases}
\sum_i \eta(e_i,e_i)^{1-g}\quad \text{ if } n=0\\
\eta(e_{i_1},e_{i_1})^{-\frac{1}{2}(2g-2+n)}\quad \text{ if } i_1=...=i_n\\
0\quad \text{ else},    
\end{cases}
\end{equation}
see \cite[Section 2.5.1]{lp}.

\section{Product targets}\label{sec:product_target}
\subsection{The Gromov-Witten CohFT}\label{sec:GW_CohFT}
For a smooth projective varieties $W$ over $\mathbb{C}$, the Gromov-Witten theory of $W$ defines a CohFT\footnote{This CohFT takes values in a suitable Novikov ring.} $\Omega^W$ with unit modelled on $H^*(W)$ equipped with the Poincar\'e pairing $\eta^W$ and $1\in H^0(W)$. The class $\Omega^W_{g,n}\in H^*(\overline{\mathcal{M}}_{g,n})\otimes H^*(W)^{\otimes n}$ is given by 
\begin{equation*}
\Omega^W_{g,n}(v_1\otimes...\otimes v_n)=\sum_{d\in H_2(W,\mathbb{Z})}Q^d \text{ft}_*\left(\prod_{i=1}^n\text{ev}_i^*(v_i)\cap [\overline{\mathcal{M}}_{g,n}(W,d)]^{vir}\right),    
\end{equation*}
where 
\begin{equation*}
\overline{\mathcal{M}}_{g,n}(W,d)
\end{equation*}
is the moduli space of genus $g$ degree $d$ $n$-pointed stable maps to $W$, $$\text{ev}_i:\overline{\mathcal{M}}_{g,n}(W,d)\to W, \quad i=1,...,n$$ are the evaluation maps and $$\text{ft}: \overline{\mathcal{M}}_{g,n}(W,d)\to \overline{\mathcal{M}}_{g,n}$$ is the forgetful map. 

Basic constructions of Gromov-Witten theory, including the virtual fundamental classes $[\overline{\mathcal{M}}_{g,n}(W,d)]^{vir}$, can be found in many references, e.g. \cite{b0}, \cite{m}.

For a class $t\in H^*(W)$, we define the {\em $t$-shifted CohFT} $\Omega^{W,t}_{g,n}$ by
\begin{equation*}
\Omega^{W,t}_{g,n}(v_1\otimes...\otimes v_n)=\sum_{m\geq 0}\frac{1}{m!}(p_m)_*\Omega^W_{g,n+m}(v_1\otimes...\otimes v_n\otimes t\otimes...\otimes t)\in H^*(\overline{\mathcal{M}}_{g,n}).
\end{equation*}
Note that the potential $\mathcal{A}^{\Omega^{W,t}}$ of the $t$-shifted Gromov-Witten CohFT of $W$ is the Gromov-Witten ancestor potential $\mathcal{A}^W_t$ of $W$, as defined in \cite[Section 5]{g2}.

\subsection{CohFT of a product}\label{sec:CohFT_product}
Suppose that the $t$-shifted Gromov-Witten CohFT $\Omega^{Y,t}$ of $Y$ is semi-simple. By Givental-Teleman classification (Theorem \ref{thm:gt}), there exists a unique $$R^Y=\text{Id}+R_1z+R_2z^2+...\in \text{End}(H^*(Y))[[z]]$$ satisfying the symplectic condition $R^Y(z)(R^Y)^{*}(-z)=\text{Id}$, such that
\begin{equation}
\Omega^{Y,t}=R^Y.\omega^{Y,t}.    
\end{equation}
Here $\omega^{Y,t}=[\Omega^{Y,t}]^0$ is the degree $0$ part of $\Omega^{Y,t}$.

Let $s\in H^*(X)$. Let $u_1,...,u_n\in H^*(X)$ and $v_1,...,v_n\in H^*(Y)$. By the product formula of Gromov-Witten invariants \cite{b}, a direct calculation show that
\begin{equation}
\Omega^{X\times Y,s\otimes 1+1\otimes t}(u_1\otimes v_1,...,u_n\otimes v_n)=\Omega^{X,s}(u_1,...,u_n)\cup \Omega^{Y,t}(v_1,...,v_n).   
\end{equation}
Therefore, we have
\begin{equation}\label{eqn:CohFT_product}
\Omega^{X\times Y,s\otimes 1+1\otimes t}=\Omega^{X,s}\Omega^{Y,t}=\Omega^{X,s}R^Y.\omega^{Y,t}=\Omega^{X,s}\sum_{\Gamma\in\mathsf{G}_{g,n}}\text{Cont}_\Gamma,
\end{equation}
where the last equality uses the description of the unit-preserving action of $R^Y$ recalled in Section \ref{sec:action}.

We examine the class $\Omega^{X,s}\text{Cont}_\Gamma$. Applying CohFT axiom (ii) to $\Omega^{X,s}$, we see that the contribution to $\Omega^{X,s}\text{Cont}_\Gamma$ from a vertex $v$ of $\Gamma$ is
\begin{equation}\label{eqn:vertex_cont}
\Omega^{X,s}_{\mathrm{g}(v),\mathrm{n}(v)}T\omega^{Y,t}_{\mathrm{g}(v),\mathrm{n}(v)}.    
\end{equation}
Since, by CohFT axiom (iii), for any $g,n,m$, we have
\begin{equation*}
\begin{split}
&\Omega^{X,s}_{g,n}(u_1\otimes...\otimes u_n)(p_m)_*\omega^{Y,t}_{g,n}(v_1\otimes...\otimes v_n\otimes T(\psi_{n+1})\otimes...\otimes T(\psi_{n+m}))\\
=&(p_m)_*(\Omega^{X,s}_{g,n+m}(u_1\otimes...\otimes u_n\otimes 1\otimes...\otimes 1)\omega^{Y,t}_{g,n}(v_1\otimes...\otimes v_n\otimes T(\psi_{n+1})\otimes...\otimes T(\psi_{n+m}))). 
\end{split}
\end{equation*}
Hence, we see that (\ref{eqn:vertex_cont}) is the translation of $\Omega^{X,s}_{\mathrm{g}(v),\mathrm{n}(v)}\omega^{Y,t}_{\mathrm{g}(v),\mathrm{n}(v)}$ by
\begin{equation*}
z\cdot 1\otimes 1-z\text{Id}\otimes (R^Y)^{-1}(1\otimes 1).   
\end{equation*}
It is clear that the contribution to $\Omega^{X,s}\text{Cont}_\Gamma$ from a leg $l$ of $\Gamma$ is 
\begin{equation}
\text{Id}\otimes (R^Y)^{-1}(\psi_l).    
\end{equation}
By CohFT axiom (ii), $\Omega^{X,s}$ contributes $(\eta^X)^{-1}$ to an edge. Therefore the contribution to $\Omega^{X,s}\text{Cont}_\Gamma$ from an edge $e$ of $\Gamma$ is
\begin{equation}
\frac{(\eta^X)^{-1}\otimes(\eta^Y)^{-1}-(\text{Id}\otimes (R^Y)^{-1})(\psi_e')((\eta^X)^{-1}\otimes(\eta^Y)^{-1})(\text{Id}\otimes (R^Y)^{-1})(\psi_e'')^{\mathbf{t}}}{\psi_e'+\psi_e''}.    
\end{equation}

Define 
\begin{equation*}
\mathsf{R}=\text{Id}\otimes R^Y\in \text{End}(H^*(X)\otimes H^*(Y))[[z]]=\text{End}(H^*(X\times Y))[[z]].
\end{equation*}
The above discussion shows that (\ref{eqn:CohFT_product}) is the unit-preserving action of $\mathsf{R}$ on the product CohFT $\Omega^{X,s}\omega^{Y,t}$.

As in (\ref{eqn:top_FT}), the topological part $\omega^{Y,t}$ can be decomposed into a sum of rank $1$ theories by using the normalized idempotent basis of $(H^*(Y),\bullet_t)$. More details about this decomposition can be found in e.g. \cite[Section 3.2.2]{gt}. Using this decomposition and the description (\ref{eqn:CohFT_product}) of $\Omega^{X\times Y, s\otimes 1+1\otimes t}$, we obtain the following formula for the Gromov-Witten ancestor potential of $X\times Y$:
\begin{equation}\label{eqn:ancestor_product}
\mathcal{A}^{X\times Y}_{s\otimes 1+1\otimes t}(\mathbf{u}\otimes \mathbf{v})=\widehat{\mathsf{R}}\prod_i\mathcal{A}^X_s(\mathbf{u}_i).    
\end{equation}
This formula, which is written in terms\footnote{It can be seen that the $R$-matrix action described in Section \ref{sec:action} corresponds to the action of the differential operator obtained by quantization of quadratic Hamiltonians in \cite{g2}.} of the quantization formulation \cite{g2}, is an extension of Givental's formula \cite[Definition 6.8]{g2} for the ancestor potential for $Y$ (i.e. $X=\text{pt}$). The potential $\mathcal{A}^{X\times Y}_{s\otimes 1+1\otimes t}$ depends on coordinates of $H^*(X)$,
\begin{equation*}
\mathbf{u}(z)=u_0+u_1z+u_2z^2+...\in H^*(X)[[z]],    
\end{equation*}
and coordinates of $H^*(Y)$,
\begin{equation*}
\mathbf{v}(z)=v_0+v_1z+v_2z^2+...\in H^*(Y)[[z]].    
\end{equation*}
The coordinates $\mathbf{u}_i(z)\in H^*(X)[[z]]$ are constructed from coordinates of $H^*(X)$ and normalized idempotents of $(H^*(Y),\bullet_t)$ in the same way as the semi-simple case, see \cite[Definition 6.8]{g2}.
By the ancestor/descendant correspondence (see \cite[Theorem 5.1]{g2}, \cite[Appendix 2]{cg}), we have the following formula for the total descendant potential of $X\times Y$:
\begin{equation}\label{eqn:descendant_potential}
\mathcal{D}^{X\times Y}=\widehat{(S^{X\times Y})^{-1}}\widehat{\mathsf{R}}\prod_i\widehat{S^{X}_s}\mathcal{D}^X,
\end{equation}
where we omit some multiplicative scalars irrelevant for our purpose.

\subsection{Virasoro constraints}\label{sec:vir}
Here we briefly review the formulation of Virasoro constriaints. Consider again a smooth projective variety $W$ over $\mathbb{C}$. Let 
\begin{equation*}
\rho^W: H^*(W)\to H^*(W)    
\end{equation*}
be the operator of multiplication by $c_1(T_W)$, the first Chern class of the tangent bundle $T_W$ of $W$. Define the {\em Hodge grading operator} of $W$, $$\mu^W: H^*(W)\to H^*(W),$$ as follows: for a homogeneous class $\phi\in H^{p,q}(W)$, define $$\mu^W(\phi)=\left(p-\frac{\text{dim}\,W}{2}\right)\phi.$$ 

The Virasoro operators for the Gromov-Witten theory of $W$ can be defined by 
\begin{equation}\label{eqn:vir_op}
L_m^W=\widehat{l_m^W}+\frac{\delta_{m,0}}{4}\text{tr}(\mu^W(\mu^W)^*), \quad m\geq -1    
\end{equation}
in terms of the quantization formulation, where
\begin{equation}
l^W_m=z^{-1/2}\left(z\frac{d}{dz}z-\mu^W z+\rho^W \right)^{m+1} z^{-1/2}.
\end{equation}
Note that for $m\geq 0$, $l_m^W=l_0^W(zl_0^W)^m$. We refer to \cite{g2} and \cite{cgt} for more details.

The Virasoro conjecture for $W$ asserts the validity of the following equalities (which we often refer to as Virasoro constraints):
\begin{equation*}
L_m^W\mathcal{D}^W=0, \quad m\geq -1.    
\end{equation*}

\subsection{Proof of Theorem \ref{thm:main}}\label{sec:pf_main}
As briefly mentioned in Section \ref{sec:intro}, we apply Givental's approach in \cite{cgt} to prove Theorem \ref{thm:main}. Fro this purpose, we must show that the operators appearing in the formula (\ref{eqn:descendant_potential}) respect gradings: namely, we want to apply the loop group covariance \cite[Proposition 1.3]{cgt}. 

Consider a product $X\times Y$ as above. We have 
\begin{equation*}
\rho^{X\times Y}=\rho^X\otimes \text{Id}+\text{Id}\otimes \rho^Y.
\end{equation*}
Also, by working with a homogeneous basis of $H^*(X\times Y)$ compatible with K\"unneth decomposition $H^*(X\times Y)=H^*(X)\otimes H^*(Y)$, we have
\begin{equation*}
\mu^{X\times Y}=\mu^X\otimes \text{Id}+\text{Id}\otimes \mu^Y.    
\end{equation*}
It follows that 
\begin{equation}\label{eqn:grading_op_XY}
l_m^{X\times Y}=l_m^X\otimes \text{Id}+\text{Id}\otimes l_m^Y.    
\end{equation}

Put $N:=\text{rank}H^*(Y)$ and $l_0^{N\text{pt}}=z\frac{d}{dz}z+\frac{1}{2}$. Then we have, by \cite[Proposition 7.7 and Theorem 8.1]{g2}, 
\begin{equation}\label{eqn:gradingY}
(S^Y)^{-1}R^Y l_0^{N\text{pt}}(R^Y)^{-1}S^Y=l_0^Y.    
\end{equation}
Together with homogeneity properties of the fundamental solutions $S^{X\times Y}$ and $S^X$ (which are consequences of the virtual dimension formula), it follows that
\begin{equation}\label{eqn:grading_XY}
(S^X\otimes \text{Id})^{-1}\mathsf{R}^{-1}S^{X\times Y}l_0^{X\times Y}(S^{X\times Y})^{-1}\mathsf{R}(S^X\otimes \text{Id})=l_0^X\otimes\text{Id}+\text{Id}\otimes l_0^{N\text{pt}}.  
\end{equation}
In other words, conjugating the operator $l_0^{X\times Y}$ across the operators on the right-hand side of (\ref{eqn:descendant_potential}) gives $$l_0^X\otimes\text{Id}+\text{Id}\otimes l_0^{N\text{pt}}=l_0^{X\times \{N \text{ points}\}}.$$ 
By the construction of Virasoro operators recalled in Section \ref{sec:vir}, we see that conjugating $L_m^{X\times Y}$ across the operators on the right-hand side of (\ref{eqn:descendant_potential}) gives 
$L_m^{X\times \{N \text{ points}\}}$. Theorem \ref{thm:main} follows.

\begin{remark}
Virasoro constraints are formulated for orbifold Gromov-Witten theory in \cite{jt}. It makes sense to ask Problem \ref{problem_product} for orbifolds. Givental-Teleman classification is for CohFTs and is applicable to orbifold Gromov-Witten theory. The product formula has been extended to orbifold Gromov-Witten theory in \cite{ajt}. Therefore, our argument for Theorem \ref{thm:main} can be extended to orbifold Gromov-Witten theory.
\end{remark}

\end{document}